\documentclass[12pt]{article}

\usepackage{fullpage,graphicx,psfrag,amsmath,amsfonts,amssymb,verbatim}
\usepackage[small,bf]{caption}
\usepackage[table]{xcolor}
\usepackage{subcaption}
\usepackage{url,hyperref}
\usepackage{listings}
\usepackage{mathtools}
\usepackage{booktabs,placeins}
\usepackage{algorithm,algpseudocode}
\usepackage{longtable}
\usepackage{color}

\newcommand{\BEAS}{\begin{eqnarray}}
\newcommand{\EEAS}{\end{eqnarray}}
\newcommand{\BEA}{\begin{eqnarray}}
\newcommand{\EEA}{\end{eqnarray}}
\newcommand{\BEQ}{\begin{equation}}
\newcommand{\EEQ}{\end{equation}}
\newcommand{\BIT}{\begin{itemize}}
\newcommand{\EIT}{\end{itemize}}
\newcommand{\BNUM}{\begin{enumerate}}
\newcommand{\ENUM}{\end{enumerate}}

\newcommand{\BA}{\begin{array}}
\newcommand{\EA}{\end{array}}

\newcommand{\eg}{{\it e.g.}}
\newcommand{\ie}{{\it i.e.}}

\newcommand{\ones}{\mathbf 1}

\newcommand{\reals}{{\mbox{\bf R}}}

\newcommand{\symm}{{\mbox{\bf S}}}  

\newcommand{\diag}{\mathop{\bf diag}}

\newcommand{\Chol}{\mathop{\bf Chol}}

\newcommand{\Expect}{\mathop{\bf E{}}}

\newcommand{\argmin}{\mathop{\rm argmin}}

\newcounter{algorithmctr}[section]
\renewcommand{\thealgorithmctr}{\thesection.\arabic{algorithmctr}}
{\refstepcounter{algorithmctr}\begin{list}{}{%
   \setlength{\rightmargin}{0.03\linewidth}%
   \setlength{\leftmargin}{0.03\linewidth}}%
   \rmfamily\small
   \item[]{\setlength{\parskip}{0ex}\hrulefill\par%
    \nopagebreak{\bfseries\textsf{Algorithm \thealgorithmctr~}}}}%
{{\setlength{\parskip}{-3ex}\nopagebreak\par\hrulefill} \end{list}}

\newcommand{\norm}[1]{\|{#1}\|}

\newcommand{\diff}{\text{d}}

\definecolor{codegreen}{rgb}{0,0.6,0}
\definecolor{codegray}{rgb}{0.5,0.5,0.5}
\definecolor{codepurple}{rgb}{0.58,0,0.82}
\definecolor{backcolour}{rgb}{0.95,0.95,0.98}

\lstdefinestyle{mystyle}{
backgroundcolor=\color{backcolour},
keywordstyle=\color{magenta},
numberstyle=\tiny\color{codegray},
stringstyle=\color{codepurple},
basicstyle=\ttfamily\small,
breakatwhitespace=false,
breaklines=true,
captionpos=t,
keepspaces=true,
numbers=left,
numbersep=5pt,
showspaces=false,
showstringspaces=false,
showtabs=false,
tabsize=2,
}

\definecolor{seagreen}{rgb}{0.18, 0.55, 0.34}
\definecolor{mediumviolet-red}{rgb}{0.78, 0.08, 0.52}
\definecolor{khaki}{rgb}{0.94, 0.9, 0.55}

\lstdefinelanguage{mypython}
{
keywords=[1]{from, import, as, assert, not, print, nonneg, PSD, axis},
keywordstyle=[1]{\color{mediumviolet-red}},
keywords=[2]{cp, lo, pl, cvxpy, Variable, Parameter,
sqrt, exp, numpy, np, Problem, Minimize, Maximize, value, solve, inner,
sum, multiply, arange, range, norm1, norm2, norm_inf, abs, square,
diagonal, outer, pos, hstack, power},
keywordstyle=[2]{\color{seagreen}},
upquote=true,
showstringspaces=false,
basicstyle=\ttfamily,
columns=fullflexible,
keepspaces=true,
emph={True,False,def,return,float,class,match,switch,len},
emphstyle={\color{seagreen}},
belowskip=1em,
aboveskip=1em,
morecomment=[l]{\#}
}

\hypersetup{ colorlinks   = true,
urlcolor     = blue,
linkcolor    = blue,
citecolor    = red
}

\graphicspath{{img/}}

\title{\LARGE \bf
Code Generation for \\Solving and Differentiating through \\Convex Optimization Problems
}
\author{Maximilian Schaller and Stephen Boyd}

\begin{document}

\maketitle
\thispagestyle{empty}
\pagestyle{empty}

\begin{abstract}
We introduce custom code generation for parametrized 
convex optimization problems that
supports evaluating the derivative of the solution with respect to the parameters,
\ie, differentiating through the optimization problem.
We extend the open source code generator CVXPYgen, which itself extends CVXPY,
a Python-embedded domain-specific language with a natural
syntax for specifying convex optimization problems,
following their mathematical description.
Our extension of CVXPYgen adds a custom C implementation to differentiate the
solution of a convex optimization problem with respect to its parameters,
together with a Python wrapper for prototyping and desktop 
(non-embedded) applications.
We give three representative application examples:
Tuning hyper-parameters in machine learning;
choosing the parameters in an approximate dynamic programming (ADP)
controller; and adjusting the parameters in an optimization based
financial trading engine via back-testing, \ie, simulation
on historical data.
While differentiating through convex optimization problems is not new,
CVXPYgen is the first tool that generates custom C code for the task,
and increases the computation speed by
about an order of magnitude in most applications, compared to CVXPYlayers,
a general-purpose tool for differentiating through convex
optimization problems.
\end{abstract}

\section{Introduction}

A convex optimization problem, parametrized by $\theta \in \Theta \subseteq \reals^d$,
can be written as
\begin{equation}\label{eq:cvx}
\begin{array}{ll}
\mbox{minimize} & f_0(x, \theta) \\
\mbox{subject to} & f_i(x, \theta) \le 0, \quad i=1,\ldots,m \\
& A(\theta)x = b(\theta),
\end{array}
\end{equation}
where $x \in \reals^n$ is the optimization variable,
$f_0$ is the objective function to be minimized, which is convex in $x$,
and $f_1,\ldots,f_m$
are inequality constraint functions that are convex in $x$~\cite{boyd2004convex}.
The parameter $\theta$ specifies data that can change, 
but is constant and given (or chosen) when we solve an instance of the problem.
We refer to the parametrized problem~\eqref{eq:cvx} as
a \emph{problem family}; when we specify a fixed value of $\theta \in \Theta$,
we refer to it as a \emph{problem instance}.
We let $x^\star$ denote an optimal point for problem~\eqref{eq:cvx}, 
assuming it exists. To emphasize its dependence on $\theta$, we write it as 
$x^\star(\theta)$.

Convex optimization is used in many domains,
including signal and image processing~\cite{mattingley2010real,zibulevsky2010l1},
machine learning~\cite{murphy2012machine,bishop2006pattern,zou2005regularization,
tibshirani1996regression,cortes1995support,hoerl1970ridge,cox1958regression},
control systems~\cite{rawlings2017model,kouvaritakis2016model,wang2015approximate,
keshavarz2014quadratic,wang2009fast,boyd1994linear,boyd1991linear,garcia1989model},
quantitative finance~\cite{palomar2025portfolio,boyd2024markowitz,boyd2017multi,narang2013inside,
lobo2007portfolio,grinold2000active,markowitz1952portfolio},
and operations research~\cite{halabian2019distributed,bertsekas2021data,bertsekas1991linear,
lefever2016convex,bertsimas2004robust}.

\subsection{Differentiating through convex optimization problems}

In many applications we are interested in the sensitivity of the solution
$x^\star$ with respect to the parameter $\theta$.
We will assume there is a unique solution for parameters near $\theta$,
and that the mapping from $\theta$ to $x^\star$ is differentiable 
with Jacobian $\partial x^\star / \partial \theta \in \reals^{n \times d}$,
evaluated at $\theta$.
We then have
\[
\Delta x^\star \approx \frac{\partial x^\star}{\partial \theta} \Delta \theta,
\]
where $\Delta \theta$ is the change in $\theta$, and
$\Delta x^\star$ is the resulting change in $x^\star$.

We make a few comments on our assumptions.
First, the solution $x^\star$ need not be unique,
and so does not define a function from $\theta$ to $x^\star$.  
Even when the solution is unique for each parameter value,
the mapping from $\theta$ to $x^\star$ need not be differentiable.
Following universal practice in machine learning, 
we simply ignore these issues.
When $x^\star$ is not unique, or when the mapping is not differentiable,
we simply use some reasonable value for the (nonexistent) derivative.
It has been observed that simple gradient (or subgradient) based methods
for optimizing parameters are tolerant of these approximations.

Approximating the change in solution with a change in parameters
can be useful by itself in some applications. As an example, 
consider a machine learning problem where we fit the parameters of 
a model to data by minimizing the sum of a convex loss function 
over the given training data.
Considering the training data as a parameter, and the solution 
as the estimated model parameters, the Jacobian above gives 
us the sensitivity of each model parameter with respect to
the training data.
In particular, these sensitivities are sometimes used to compute a risk
estimate for the learned model parameters
\cite{nobel2023tractable, nobel2024randalo}.
As another example, suppose we model some economic variables 
(\eg, consumption, demand for products, trades) as maximizing
a concave utility function that depends on parameters.
The Jacobian here directly gives us an approximation of 
the change in demand (say) when the utility parameters change
\cite{wainwright2005fundamental}.

\pagebreak
\subsection{Autodifferentiation framework}
The solution map derivative is much more useful when it is part
of an autodifferentiation system such as JAX~\cite{jax2018github},
PyTorch~\cite{paszke2017automatic}, or Tensorflow~\cite{abadi2016tensorflow}.
We consider a scalar function that is described by its compution graph,
which can include standard operations and functions, as well
the solution of one or more convex optimization problems.
We can compute the gradient of this function automatically,
and this can be used for applications such as tuning or optimizing
the performance of a system.
We give a few simple generic examples here.

In machine learning, we fit model parameters (also called weights or
coefficents) using convex
optimization, but we may have other hyper-parameters (\eg, that scale 
regularization terms) that we would like to tune to get good 
performance on an unseen, out-of-sample test data set~\cite{murphy2012machine}.
The scalar function that we differentiate is the loss function computed with
test data, and we differentiate with respect to the hyper-parameters
of the machine learning model.
A similar situation occurs in finance, where the actual trades to execute
are determined by solving a parameterized problem~\cite{boyd2017multi}, which also
contains a number of hyper-parameters that set limits on the portfolio or 
trading, or scale objective terms, and our goal is to obtain good 
performance on a simulation that uses historical data, \ie, a back-test.
In this case, the scalar function that we optimize might be a metric like
the realized portfolio return or the Sharpe ratio~\cite{ledoit2008robust},
and we differentiate with respect to
the hyper-parameters of the portfolio construction model.

\subsection{Related work}

\paragraph{Differentiating through convex optimization problems.}
There are two main classes of methods to differentiate the mapping from
problem parameters to the solution.
First, autodifferentiation software like
PyTorch~\cite{paszke2017automatic} and Tensorflow~\cite{abadi2016tensorflow},
as commonly used in backpropagation for
machine learning, would differentiate through all instructions of an iterative
optimization algorithm. While these tools are commonly used in the
deep learning domain, they neglect the structure and optimality conditions of,
particularly, convex optimization problems.

Second, as less of a brute-force approach,
one can directly differentiate through the optimality conditions of a convex
optimization problem, also referred to as
\emph{argmin differentiation}~\cite{agrawal2019differentiating,amos2017optnet}.
CVXPYlayers~\cite{agrawal2019differentiable} is based on prior work on differentiating
through a cone program, called diffcp~\cite{agrawal2019differentiating},
and provides interfaces to PyTorch and Tensorflow.
OptNet~\cite{amos2017optnet} addresses quadratic programs.
Least squares auto-tuning~\cite{barratt2021least} is a specialized hyper-parameter
tuning framework for least squares problems, \ie, a subclass of quadratic programs,
with parametrized problem data.
It admits flexible tuning objectives and tuning regularizers (called hyper-hyper-parameters),
and computes their gradient with respect to the parameters of the least squares problem.
PyEPO~\cite{tang2024pyepo} combines various autodifferentiation and argmin
differentiation tools into one software suite.

\paragraph{Tuning systems that involve convex optimization.}
Differentiating through a convex problem is useful in a broad array
of applications sometimes called predict-then-optimize. In these applications
we make some forecast or prediction (using convex optimization or other
machine learning models),
and then take some action (using convex optimization),
and we are interested in the gradient of parameters appearing both in the 
predictor and the action policy~\cite{tang2024pyepo,elmachtoub2022smart}.
In some cases, the prediction and optimization steps are fused
into one black-box model that maps features to actions,
called \emph{learning to optimize from features}~\cite{kotary2024learning},
which avoids differentiating through an optimization problem.
This approximation is not necessary when it is possible to differentiate
through the optimization problem in an easy and fast way,
which motivates this work.

We also mention that when the dimension of $\theta$ is small enough,
the parameters can be optimized
using zero-order or derivative-free methods, which do not require the gradient
of the overall metric with respect to the
parameters~\cite{agrawal2020learning,agrawal2019differentiable}.
Examples include Optuna~\cite{akiba2019optuna},
HOLA~\cite{maher2022light}, for tuning parameters with respect to some 
overall metric, and Hyperband~\cite{li2018hyperband}, which
dynamically allocates resources for efficient hyper-parameter search in 
machine learning.
Even when a tuning problem can be reasonably carried out using
derivative-free methods, the ability to evaluate the gradient can 
give faster convergence with fewer evaluations.

\paragraph{Domain-specific languages for convex optimization.}

Argmin differentiation tools like CVXPYlayers admit convex optimization problems
that are specified in a domain-specific language (DSL).
Such systems allow the user to specify the functions $f_i$ and $A$ and $b$,
in a simple format that closely follows the mathematical description of the problem.
Examples include YALMIP~\cite{lofberg2004yalmip} and
CVX~\cite{grant2014cvx} (in Matlab), CVXPY~\cite{diamond2016cvxpy} (in Python),
which CVXPYlayers is based on,
Convex.jl~\cite{convexjl} and JuMP~\cite{Dunning2017jump} (in Julia),
and CVXR~\cite{fu2017cvxr} (in R).
We focus on CVXPY, which also supports the declaration of parameters,
enabling it to specify problem families, not just problem instances.

DSLs parse the problem description and translate (canonicalize)
it to an equivalent problem
that is suitable for a solver that handles some generic class of problems,
such as linear programs (LPs), quadratic programs (QPs), second-order cone programs (SOCPs),
semidefinite programs (SDPs), and others such as exponential cone programs~\cite{boyd2004convex}.
Our work focuses on LPs and QPs.
After the canonicalized problem is solved, a solution of the original problem is retrieved
from a solution of the canonicalized problem.

It is useful to think of the whole process as a function that maps $\theta$,
the parameter that specifies the problem instance, into $x^\star$, an optimal
value of the variable. With a DSL, this process consists of three steps.
First the original problem description is canonicalized to a problem in some 
standard (canonical) form;
then the canonicalized problem is solved using a solver; and finally, a solution
of the original problem is retrieved from a solution of the canonicalized problem.
When differentiating through the problem,
the sequence of canonicalization, canonical solving, and retrieval is reversed.
Reverse retrieval is followed by canonical differentiating and reverse canonicalization.

\begin{figure}
\centering
\begin{subfigure}{\columnwidth}
\centering
\includegraphics[width=0.8\linewidth]{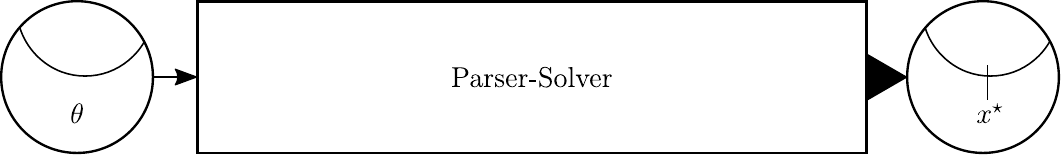}
\caption{Parser-solver calculating solution $x^\star$ for 
problem instance with parameter $\theta$.}
\label{fig:flow_general}
\end{subfigure}
\par\vspace{10pt}
\begin{subfigure}{\columnwidth}
\centering
\includegraphics[width=0.8\linewidth]{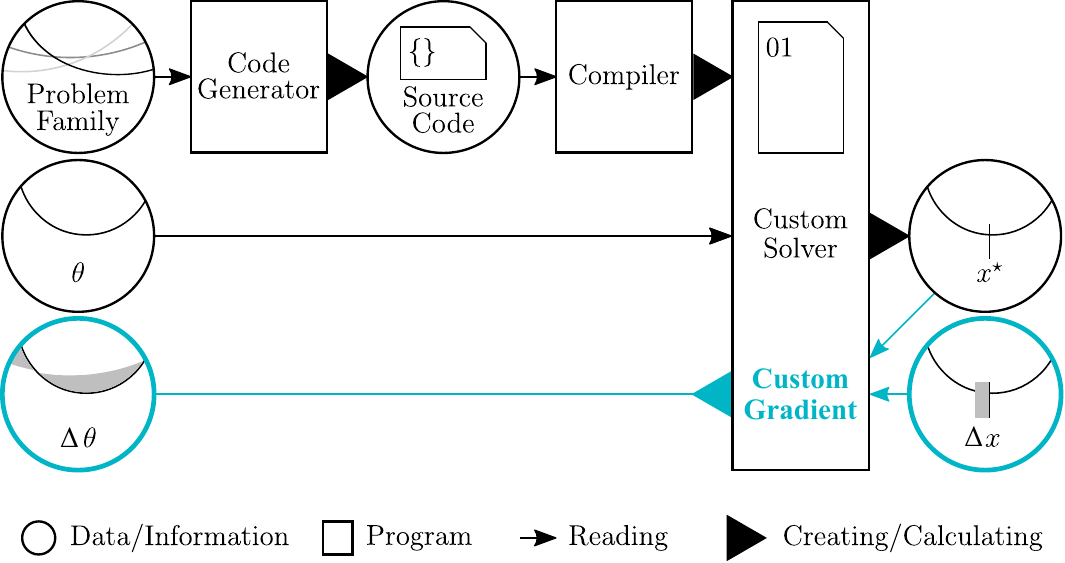}
\caption{Source code generation for problem family, followed by compilation to 
custom solver and custom gradient computation (new, signified with blue color).
The compiled solver computes a solution $x^\star$ to the problem instance with
parameter $\theta$.
The compiled differentiator computes the gradient $\Delta \theta$
given $\Delta x$.}
\label{fig:flow_custom}
\end{subfigure}
\caption{Comparison of convex optimization problem parsing and
solving/differentiating approaches.}
\vspace{-1em}
\label{fig:flow}
\end{figure}

Most DSLs are organized as
\emph{parser-solvers}, which carry out the canonicalization each time the problem is
solved (with different parameter values).
This simple setting is illustrated in figure~\ref{fig:flow_general}.

\paragraph{Code generation for convex optimization.}

We are interested in applications where we solve many instances of the
problem, possibly in an embedded application with hard real-time constraints, or a
non-embedded application with limited compute.
For such applications, a \emph{code generator} makes more sense.

A code generator takes as input a description of a problem family, and generates
specialized source code for that specific family.
That source code is then compiled, and we have an efficient solver for
the specific family. In this work, we add a program that efficiently computes
the gradient of the parameter-solution mapping.
The overall workflow is illustrated in figure~\ref{fig:flow_custom}.
The compiled solver and differentiator have a number of advantages over parser-solvers.
First, by caching canonicalization and exploiting the problem structure,
the compiled solver and the compiled differentiator are faster.
Second, the compiled solver and in some applications also the
compiled differentiator can be deployed in embedded systems,
fulfilling rules for safety-critical code~\cite{holzmann2006power}.

\FloatBarrier

\subsection{Contribution}\label{sec:contribution}

In this paper, we extend the code generator CVXPYgen~\cite{schaller2022embedded}
to produce source code for differentiating
the parameter-solution mapping of convex optimization problems
that can be reduced to QPs.
We allow for the use of any canonical solver that is supported by CVXPYgen,
including conic solvers.
We combine existing theory on
differentiating through the optimality conditions of a QP~\cite{amos2017optnet} with
low-rank updates to the factorization of quasidefinite
systems~\cite{davis1999modifying, davis2005row}
to enable very fast repeated differentiation.
Along with the generated C code, we compile two Python interfaces, one
for use with CVXPY and one for use with CVXPYlayers.
To the best of our knowledge, CVXPYgen is the first code generator for
convex optimization that supports differentiation.

We give three examples, tuning the hyper-parameters and
feature engineering parameters of a machine learning model,
tuning the controller weights of an approximate dynamic programming controller,
and adjusting the parameters in a financial trading engine.
CVXPYgen accelerates these applications by around an order of magnitude.

\subsection{Outline}

The remainder of this paper is structured as follows.
In \S\ref{sec:CVXPYgen} we describe, at a high level, how CVXPYgen generates
code to differentiate through convex optimization problems.
In \S\ref{sec:tuning_framework} we describe the generic system tuning framework that
runs the CVXPYgen solvers and differentiators, and in
\S\ref{sec:examples}, we illustrate how we use it for three realistic examples
that compare the performance of CVXPYgen to CVXPYlayers.
We conclude the paper in \S\ref{sec:conclusion}.

\section{Differentiating with CVXPYgen}\label{sec:CVXPYgen}

CVXPYgen is an open-source code generator, based on the
Python-embedded domain-specific language CVXPY.
While CVXPY treats all typical conic programs and CVXPYgen generates code
to solve LPs, QPs, and SOCPs, we focus on differentiating through problems
that can be reduced to a QP, \ie, LPs and QPs.

\subsection{Disciplined parametrized programming}
The CVXPY language uses disciplined convex programming (DCP) to allow
for modeling instructions that are very close to the mathematical problem description
and to verify convexity in a systematic way~\cite{diamond2016cvxpy}.
Disciplined parametrized programming (DPP) is an extension to the DCP rules
for modeling convex optimization problems.
While a DCP problem is readily canonicalized, a DPP problem is readily canonicalized
with \emph{affine} mappings from the user-defined parameters to the canonical parameters.
Similarly, the mapping from a canonical solution back to a solution to the user-defined problem
is affine for DPP problems~\cite{agrawal2019differentiable}.
DPP imposes mild restrictions on how parameters enter the problem
expressions. In short, if all parameters enter the problem expressions in
an affine way, the problem is DPP.
We model all example problems in \S\ref{sec:examples} DPP and illustrate standard
modifications to make DCP problems DPP.
Details on the DCP and DPP rules can be found at \url{https://www.cvxpy.org}.

\subsection{Differentiating through parametrized problems}
Differentiating through a DPP problem consists of three steps:
affine parameter canonicalization, canonical solving,
and affine solution retrieval,
\[
\tilde{\theta} = C  \theta + c, \qquad
\tilde x^\star = \mathcal{S} (\tilde \theta), \qquad
x^\star = R \tilde x^\star + r,
\]
where $\theta$ is the user-defined parameter, $C$ and $R$ are sparse matrices,
and $\mathcal{S}(\cdot)$ is the canoncial solver. We mark the canonical
parameter and solution with a tilde.

In this work, we want to propagate a gradient in terms of
the current solution, called $\Delta x$,
to a gradient in the parameters $\Delta \theta$.
This mapping is symmetric to the solution
mapping~\cite{agrawal2019differentiable}, \ie,
\[
\Delta \tilde{x} = R^T \Delta x, \qquad
\Delta \tilde{\theta} = (\mathsf{D}^T \mathcal{S})
(\Delta \tilde x; \tilde x^\star, \tilde \theta), \qquad
\Delta \theta = C^T \Delta  \tilde{\theta},
\]
and we can simply re-use the descriptions of $R$ and $C$
that CVXPYgen has already extracted for solving the problem.
The following section explains how we (re)compute the canonical derivative
$(\mathsf{D}^T \mathcal{S})(\Delta \tilde x; \tilde x^\star, \tilde \theta)$ efficiently.

\subsection{Differentiating through canonical solver}

We focus on differentiating through LPs and QPs,
\ie, problems that can be reduced to the QP standard form
\[
\begin{array}{ll}
\mbox{minimize} & (1/2) \tilde x^T P \tilde x + q^T \tilde x \\
\mbox{subject to} & l \le A \tilde x \le u,
\end{array}
\]
as used by the OSQP solver~\cite{stellato2020osqp}.
The variable is $\tilde x \in \reals^{\tilde n}$ and all other symbols are
parameters.
The objective is parametrized with
$P \in \symm^{\tilde n}_{+}$, where $\symm^{\tilde n}_{+}$ is the set of
symmetric positive semidefinite matrices, and $q \in \reals^{\tilde n}$.
The constraints are parametrized with $A \in \reals^{m \times {\tilde n}}$,
$l \in \reals^m \cup \{-\infty\}$, and $u \in \reals^m \cup \{\infty\}$.
If an entry of $l$ or $u$ is $-\infty$ or $\infty$, respectively,
it means there is no constraint.
A pair of equal entries $l_i = u_i$ represents an equality constraint.

We closely follow the approach that is used in
OptNet~\cite{amos2017optnet}.
We denote by $A_\mathcal{C}$ the row slice of $A$
that contains all rows $A_i$
for which the lower or upper constraint is active
at optimality, \ie, it holds that $A_i\tilde x = l_i$ or
$A_i\tilde x = u_i$ (or both, in the case of an equality constraint).
We omit the superscript $\star$ for brevity.
We set $b_\mathcal{C}$ to contain
the entries of $l$ or $u$ at the active constraints indices, \ie,
$A_\mathcal{C} \tilde x = b_\mathcal{C}$.
Then, the solution is characterized by the KKT system
\BEQ\label{eq:qp_kkt}
\begin{bmatrix}
P & A_\mathcal{C}^T \\
A_\mathcal{C} & 0
\end{bmatrix}
\begin{bmatrix}
\tilde x \\
\tilde y_\mathcal{C}
\end{bmatrix}
=
\begin{bmatrix}
-q \\
b_\mathcal{C}
\end{bmatrix},
\EEQ
where $\tilde y_\mathcal{C}$ is the slice of the
dual variable corresponding to the active constraints.
Note that we use the sign of $\tilde y$
to determine constraint activity.
(The first block row of system~\eqref{eq:qp_kkt}
corresponds to stationarity of the Lagrangian,
and the second block row corresponds to primal feasibility.)

We take the differential of \eqref{eq:qp_kkt} and re-group the terms as
\[
\begin{bmatrix}
P & A_\mathcal{C}^T \\
A_\mathcal{C} & 0
\end{bmatrix}
\begin{bmatrix}
\diff \tilde x \\
\diff \tilde y_\mathcal{C}
\end{bmatrix}
=
\begin{bmatrix}
-\diff P \tilde x - \diff A_\mathcal{C}^T \tilde y_\mathcal{C} -\diff q \\
-\diff A_\mathcal{C} \tilde x  + \diff b_\mathcal{C}
\end{bmatrix}.
\]
We introduce
\BEQ\label{eq:dsystem}
\begin{bmatrix}
d_x \\
d_y
\end{bmatrix}
=
\begin{bmatrix}
P & A_\mathcal{C}^T \\
A_\mathcal{C} & 0
\end{bmatrix}^{-1}
\begin{bmatrix}
\Delta \tilde x \\
0
\end{bmatrix},
\EEQ
where the righthand side is evaluated using algorithm~\ref{alg:ir}.
To avoid singularity of the linear system, we regularize the matrix diagonal
with a small $\epsilon > 0$, solve the regularized system and add
$N^\text{refine}$ steps of iterative refinement~\cite{carson2018accelerating, higham1997iterative}
to correct for the effect of
$\epsilon$ on the solution.
\begin{algorithm}
\caption{Regularized system solve}
\label{alg:ir}
\begin{algorithmic}[1]
\State Initialize $P$, $A_\mathcal{C}$, $\Delta \tilde x$
\State
$K_\mathcal{C} = \begin{bmatrix}
P & A_\mathcal{C}^T \\
A_\mathcal{C} & 0
\end{bmatrix}$,
\quad
$K_\mathcal{C}^\epsilon = K_\mathcal{C} +
\begin{bmatrix} \epsilon I & 0 \\ 0 & -\epsilon I \end{bmatrix}$,
\quad
$r = \begin{bmatrix} \Delta \tilde x \\ 0 \end{bmatrix}$
\State $z = (K_\mathcal{C}^\epsilon)^{-1} r$
\For{$i = 1$ to $N^\text{refine}$}
\State $\delta_r = r - K_\mathcal{C} z$
\State $\delta_z = (K_\mathcal{C}^\epsilon)^{-1} \delta_r$
\State $z \leftarrow z + \delta_z$
\EndFor
\State $(d_x, d_y) = z$
\end{algorithmic}
\end{algorithm}
\FloatBarrier
The regularization strength $\epsilon = 10^{-6}$ and
$N^\text{refine} = 3$ iterations of iterative refinement work well in most practical cases.

Ultimately, the gradients in the QP parameters are
\[
\Delta P = -(1/2)  ( d_x \tilde x^T + \tilde x d_x^T), \qquad
\Delta q = -d_x,  \qquad
\Delta A_\mathcal{C} = - (d_y \tilde x^T + \tilde y_\mathcal{C} d_x^T),
\qquad \Delta b_\mathcal{C} = d_y.
\]
We copy the rows of $\Delta A_\mathcal{C}$ and $\Delta b_\mathcal{C}$
into the corresponding rows of $\Delta A$ and $\Delta b$,
respectively, and set all other entries of $\Delta A$ and $\Delta b$ to zero.

\paragraph{Low-rank updates to factorization of linear system.}
For the quasidefinite matrix $K_\mathcal{C}^\epsilon$ used in algorithm~\ref{alg:ir},
there always exists an LDL-factorization~\cite{vanderbei1995symmetric},
\[
L_\mathcal{C}D_\mathcal{C}L_\mathcal{C}^T =
K_\mathcal{C}^\epsilon = 
\begin{bmatrix}
P + \epsilon I & A_\mathcal{C}^T \\
A_\mathcal{C} & -\epsilon I
\end{bmatrix}.
\]
If the entries of $P$ or $A$ change, we perform a full re-factorization.
Otherwise, we use the fact that the factors $L_\mathcal{C}$ and $D_\mathcal{C}$
change with the set of active constraints, denoted by $\mathcal{C}$.
For the constraints that switch
from inactive to active or vice-versa, we perform a sequence of rank-1 updates
to $L_\mathcal{C}$ and $D_\mathcal{C}$.

We re-write the LDL-factorization as
\[
\begin{bmatrix}
L_{11} & 0 & 0 \\
\bar{l}_{12}^T & 1 & 0 \\
L_{31} & \bar{l}_{32} & \bar{L}_{33}
\end{bmatrix}
\begin{bmatrix}
D_{11} & 0 & 0 \\
0 & \bar{d}_{22} & 0 \\
0 & 0 & \bar{D}_{33}
\end{bmatrix}
\begin{bmatrix}
L_{11}^T & \bar{l}_{12} & L_{31}^T \\
0 & 1 & \bar{l}_{32}^T \\
0 & 0 & \bar{L}_{33}^T
\end{bmatrix}
=
\begin{bmatrix}
K_{11} & \bar{k}_{12} & K_{31}^T \\
\bar{k}_{12}^T & \bar{k}_{22} & \bar{k}_{32}^T \\
K_{31} & \bar{k}_{32} & K_{33}
\end{bmatrix},
\]
where lowercase symbols marked with a bar denote
row/column combinations that are added or deleted.
Uppercase symbols marked with a bar are
sub-matrices that will be altered due to the addition or deletion.

For a constraint that switches from inactive to active, we add the respective row/column
combination to $K_\mathcal{C}^\epsilon$ and run algorithm
\ref{alg:addition}.
\begin{algorithm}
\caption{Row/column addition (variant of algorithm 1 in~\cite{davis2005row})}
\label{alg:addition}
\begin{algorithmic}[1]
\State Solve the lower triangular system
$L_{11} D_{11} \bar{l}_{12} = \bar{k}_{12}$ for $\bar{l}_{12}$
\State $\bar{d}_{22} = \bar{k}_{22} - \bar{l}_{12}^T D_{11} \bar{l}_{12}$
\State $\bar{l}_{32} = (\bar{k}_{32} - L_{31} D_{11} \bar{l}_{12}) / \bar{d}_{22}$
\State $w = \bar{l}_{32} (-\bar{d}_{22})^{1/2}$
\State Perform the rank-1 downdate
$\bar{L}_{33} \bar{D}_{33} \bar{L}_{33}^T = L_{33} D_{33} L_{33}^T - w w^T$
according to algorithm 5 in~\cite{davis1999modifying}
\end{algorithmic}
\end{algorithm}
If a constraint switches from active to inactive, we run algorithm
\ref{alg:deletion} for row/column deletion.
\begin{algorithm}
\caption{Row/column deletion (variant of algorithm 2 in~\cite{davis2005row})}
\label{alg:deletion}
\begin{algorithmic}[1]
\State $w = \bar{l}_{32} (-\bar{d}_{22})^{1/2}$
\State $\bar{l}_{12} = 0$
\State $\bar{d}_{22} = 1$
\State $\bar{l}_{32} =0$
\State Perform the rank-1 update
$\bar{L}_{33} \bar{D}_{33} \bar{L}_{33}^T = L_{33} D_{33} L_{33}^T + w w^T$
according to algorithm 5 in~\cite{davis1999modifying}
\end{algorithmic}
\end{algorithm}
All steps of the row/column addition and deletion algorithms operate on the
sparse matrices $K_\mathcal{C}^\epsilon$ and $L_\mathcal{C}$ stored in
compressed sparse column format. The diagonal matrix
$D_\mathcal{C}$ is stored as an array of diagonal entries.
Note that $(-\bar{d}_{22})^{1/2}$ is always real because we run
algorithms~\ref{alg:addition} and~\ref{alg:deletion} only for row/column combinations
that are in the lower and right parts of
$K_\mathcal{C}^\epsilon$ (where $A_\mathcal{C}$ changes),
for which the diagonal entries of $D_\mathcal{C}$ are all negative by
quasidefiniteness of $K_\mathcal{C}^\epsilon$.

It is important to note that all steps, including step 5 in both algorithms,
are of at most quadratic complexity, whereas a full re-factorization
would be of cubic complexity.
When only a few constraints switch their activity, this procedure is
considerably faster than full re-factorizations.
This is demonstrated in \S\ref{sec:examples} for three practical cases.
In the worst case where all constraints switch to active or inactive
between two consecutive solves, the complexity returns to cubic.

Open source code and full documentation for CVXPYgen and its 
differentiation feature is available at
\begin{center}
\url{https://pypi.org/project/cvxpygen}.
\end{center}

\section{System tuning framework}\label{sec:tuning_framework}
We present a generic tuning method for systems of the form
\[
p = \Gamma(\omega)
\]
where $p \in \reals$ is a performance objective, $\Gamma$ evaluates the system,
which includes many solves of the convex optimization problem~\eqref{eq:cvx},
possibly sequentially, and $\omega \in \Omega \subseteq \reals^p$
is a \emph{design}, where $\Omega$ is the design space, \ie, the set
of admissible designs.
Then, we describe in detail what $\Gamma$ and $\omega$ are, for
two important classes of system tuning.

\subsection{A generic tuning method}\label{sec:pgd}

We compute the gradient $\nabla \Gamma(\omega)$ using the chain rule
and our ability to differentiate through DPP problems.
We use $\nabla \Gamma(\omega)$ in a simple
projected gradient method~\cite{bertsekas1997nonlinear,calamai1987projected}
to optimize the design $\omega$.

We use Euclidean projections onto the design space $\Omega$,
denoted by $\Pi(\cdot)$, and
a simple line search to guarantee that the algorithm is a descent method.
If the performance is improved with the current step size, we use it for the current iteration
and increase it by a constant factor $\beta > 1$ for the next iteration.
Otherwise, we repeatedly shrink the step size by a constant factor $\eta > 1$
until the performance is improved.
The simple generic design method we use is given in algorithm~\ref{alg:pgd}.

\begin{algorithm}
\caption{Projected gradient descent}\label{alg:pgd}
\begin{algorithmic}[1]
\State Initialize $\omega^0$, $\alpha^0$, $k=0$
\Repeat
\State $\hat \omega = \Pi (\omega^k - \alpha^k \nabla \Gamma(\omega^k))$
\Comment tentative update
\If{$\Gamma(\hat \omega) < \Gamma(\omega^k)$}
\State $\omega^{k+1} = \hat \omega$, $\alpha^{k+1} = \beta \alpha^k$
\Comment accept update and increase step size
\Else
\State $\alpha^k \gets \alpha^k / \eta$, go to step 3
\Comment shrink step size and re-evaluate
\EndIf
\State $k \leftarrow k + 1$
\Until{$\norm{\omega^k -
\Pi (\omega^k - \alpha^k \nabla \Gamma(\omega^k))}_2 \leq
\epsilon^\text{rel} \norm{\omega^k}_2 + \epsilon^\text{abs}$}
\end{algorithmic}
\end{algorithm}

\FloatBarrier

Note that algorithm \ref{alg:pgd} assumes that $\Gamma(\omega)$ is to be minimized.
If $\Gamma(\omega)$ is to be maximized, replace it with $-\Gamma(\omega)$.

\paragraph{Initialization.}
We initialize $\omega^0$ to a value that is typical for the
respective application and $\alpha^0$ with the clipped Polyak step size
\[
\alpha^0 = \min \left\{
\frac{\Gamma(\omega^0) - \hat p}{\norm{\nabla \Gamma(\omega^0)}_2^2},
1
\right\},
\]
where $\hat p$ is an estimate for the optimal value of the performance objective.
We clip the step size at 1 to avoid too large initial steps due to local concavity.
The algorithm is not particularly dependent on the line search parameters $\beta$
and $\eta$. Reasonable choices are, \eg, $\beta = 1.2$ and $\eta = 1.5$.

\paragraph{Stopping criterion.}
We stop the algorithm as soon as the termination criterion
\[
\norm{\omega^k -
	\Pi (\omega^k - \alpha^k \nabla \Gamma(\omega^k))}_2 \leq
\epsilon^\text{rel} \norm{\omega^k}_2 + \epsilon^\text{abs}
\]
with $\epsilon^\text{rel}, \epsilon^\text{abs} > 0$ is met.
This is also referred to as the \emph{projected gradient} being small.
When $\Gamma$ is convex, this corresponds to the first-order optimality condition, \ie,
the gradient $\nabla \Gamma(\omega^k)$ lying in (or close to) the normal cone
to $\Omega$ at $\omega^k$~\cite{boyd2004convex}.
Depending on the application, the stopping tolerances $\epsilon^\text{rel}$
and $\epsilon^\text{abs}$ might range between $10^{-2}$ and $10^{-6}$.

\subsection{Tuning hyper-parameters of machine learning models}\label{sec:mltuning}
We call the data points used for \emph{training} a machine learning model
$(z_1, y_1), \ldots, (z_N, y_N) \in \mathcal{D}$,
where each data point consists of features $z_i$ and output $y_i$.
For the design $\omega$ of the machine learning model, we consider any hyper-parameters,
including pre-processing parameters that determine how the data
$(z_1, y_1), \ldots, (z_N, y_N)$ is modified before fitting the model.

For any choice of $\omega$,
we find the model weights $\beta \in \reals^n$ as
\[
\beta^\star (\omega) = \argmin_\beta
\frac{1}{N} \sum_{i=1}^N \ell(z_i, y_i, \beta, \omega) + r(\beta, \omega),
\]
where $\ell : \mathcal{D} \times \reals^n \times \Omega \rightarrow \reals$
is the training loss
function and $r :  \reals^n \times \Omega \rightarrow \reals$ is the regularizer.
While the entries of $\beta$ are oftentimes referred to as \emph{model parameters} in the
machine learning literature, we call them \emph{weights} to make clear that
they enter the above optimization problem as variables (and not as parameters of
the optimization problem).
Both $\ell$ and $r$ are parametrized by the design $\omega$.
In the case of $\ell$,
the design $\omega$ might enter in terms of pre-processing parameters
like thresholds for processing outliers.
In the case of $r$, the design $\omega$ might enter
as the scaling of the regularization term.
In the remainder of this work, we consider $\ell$ and $r$ that are convex and quadratic,
admitting the differentiation method described in \S\ref{sec:CVXPYgen}.

We choose $\omega$ to minimize the validation loss
\[
p = \Gamma(\omega) = \frac{1}{N^\text{val}} \sum_{i=1}^{N^\text{val}}
\ell^\text{val} (z^\text{val}_i, y^\text{val}_i, \beta^\star(\omega), \omega),
\]
where $\ell^\text{val} : \mathcal{D} \times \reals^n \times \Omega \rightarrow \reals$
is the validation loss function, with validation data
$(z^\text{val}_i, y^\text{val}_i)$
that is different from and ideally uncorrelated with the training data.
Here, $\ell^\text{val}$ need not be convex (or quadratic),
since we use the projected gradient method described in \S\ref{sec:pgd}.
Note that the design $\omega$ enters $\ell^\text{val}$ both directly
(for example as pre-processing parameters) and
through the optimal model parameters $\beta^\star(\omega)$.

If $\ell^\text{val}$ is the squared error (between data and model output),
then the alternative performance objective
\[
\bar p = \bar\Gamma(\omega) = \left(\frac{1}{N^\text{val}} \sum_{i=1}^{N^\text{val}}
\ell^\text{val} (z^\text{val}_i, y^\text{val}_i, \beta^\star(\omega), \omega)\right)^{1/2},
\]
is more meaningful, as it resembles the root mean square error (RMSE).

\paragraph{Cross validation.}
For better generalization of the
optimized $\omega$ to unseen data, we can employ
cross validation (CV)~\cite{shao1993linear, hastie2009elements}.
We split the set of data points into $J$ partitions or \emph{folds}
(typically equally sized) and train the model $J$ times.
Every time, we take $J-1$ folds as training data and $1$ fold as validation data.
We compute the performance $p_j$ as described above,
where the subscript $j$ denotes that
the validation data $(z^\text{val}_i, y^\text{val}_i)$ is that of the $j$th fold.
Then, we average these over all folds as
\[
p^\text{CV} = \frac{1}{J} \sum_{j=1}^J p_j(\omega).
\]
This is usually referred to as the cross validation loss.
Similarly, if $\ell^\text{val}$ is the squared error, we compute the cross-validated
RMSE $\bar p^\text{CV}$ by averaging all $\bar p_j$.

\subsection{Tuning the weights of convex optimization control policies}
\label{sec:cocptuning}

We consider a convex optimization control policy (COCP) that
determines a control input $u \in \mathcal{U} \subseteq \reals^m$ that is applied
to a dynamical system with state $x \in \mathcal{X} \subseteq \reals^n$,
by solving the convex optimization problem
\[
u = \phi(x; \omega) = \argmin_{u \in \mathcal{U}} \ell(x, u, \omega).
\]
Here, $\phi : \mathcal{X} \times \Omega \rightarrow \mathcal{U}$ is a family of
control policies, parametrized by $\omega$,
and $x$ is the current measurement (or estimate) of the state of the dynamical system.
The loss function $\ell :\mathcal{X} \times \mathcal{U} \times \Omega \rightarrow \reals$
is parametrized by the design $\omega$,
which might involve controller weights, for example.

We choose $\omega$ to minimize the closed-loop loss
\[
p = \Gamma(\omega) = \ell^\text{cl} (x_0, \omega),
\]
where the loss function
$\ell^\text{cl} : \mathcal{X} \times \Omega \rightarrow \reals$ involves a
simulation or real experiment starting from the initial state $x_0$ and using
the control policy $u = \phi(x; \omega)$.

\section{Numerical experiments}\label{sec:examples}

In this section we present three numerical examples, 
comparing the solve and differentiation speed of CVXPYgen with CVXPYlayers,
for system tuning with the framework presented in \S\ref{sec:pgd}.
In all three cases we take OSQP as the canonical solver for CVXPYgen
and Clarabel as the canonical solver for CVXPYlayers (since it only supports
conic solvers).
The code that was used for the experiments is available at
\begin{center}
\url{https://github.com/cvxgrp/cvxpygen}.
\end{center}

\subsection{Elastic net regression with winsorized features}
\label{sec:example_elastic}

We consider a linear regression model with elastic net
regularization~\cite{zou2005regularization},
which is a sum of ridge (sum squares)~\cite{hoerl1970ridge} and
lasso (sum absolute)~\cite{tibshirani1996regression} regularization, each
of which has a scaling parameter. In addition, we clip or \emph{winsorize} each feature
at some specified level to mitigate the problem of feature outliers.
The clipping levels for each feature are also parameters.

We consider $m$ data points with one observation and $n$ features each.
We start with a set of observations $y_i \in \reals$ and
raw features $z_i \in \reals^n$,
where $z_{i,j}$ is the $j$th feature for the $i$th observation, subject to outliers.
We obtain the winsorized features $x_i \in \reals^n$ by
clipping the $n$ components at winsorization levels $w \in \reals^n_{++}$
(part of the design) as
\BEQ\label{eq:winsorize}
x_{i,j} (z_{i,j}; w) = \min \{ \max \{ z_{i,j}, -w_j \}, w_j \},
\quad i = 1,\ldots,m, \quad j = 1,\ldots,n.
\EEQ

The training loss function and regularizer as in \S\ref{sec:mltuning} are
\[
\ell(z_i, y_i, \beta, \omega) = (x_i(z_i; w)^T \beta - y_i)^2, \quad
r(\beta, \omega) = \lambda \norm{\beta}_2^2 + \gamma \norm{\beta}_1,
\]
where $\lambda \geq 0$ and $\gamma \geq 0$ are the the ridge and
lasso regularization factors, respectively.
Since the model performance depends mainly on the orders of magnitude
of $\lambda$ and $\gamma$, we write them as
$\lambda = 10^\mu$ and $\gamma = 10^\nu$.

Together with the winsorization levels $w$, the design vector
becomes $\omega = (w, \mu, \nu) \in \reals^n_{++} \times \reals^2$.
We allow to clip $z_i$ between 1 and 3 standard deviations and
search the elastic net weights across 7 orders of magnitude.
We assume that the entries of $z_i$ are approximately centered and scaled, \ie,
have zero mean and standard deviation 1.
The design space becomes
\[
\Omega = [1, 3]^n \times [-3, 3]^2.
\]

The validation loss function $\ell^\text{val}$ is identical to the training loss
function $\ell$ and the performance objective is the
cross-validated RMSE $\bar p^\text{CV}$ from \S\ref{sec:mltuning}.

\paragraph{Code generation.}

Figure~\ref{code:elastic_CVXPYgen} shows how to generate code for this problem.
\begin{figure}
\lstset{language=Python,
numbers=left,
xleftmargin=0.07\columnwidth,
linewidth=\columnwidth}
\begin{lstlisting}[frame=lines]
import cvxpy as cp
from cvxpygen import cpg

# model problem
beta = cp.Variable(n, name='beta')
X = cp.Parameter((m-m//J, n), name='X')
l = cp.Parameter(nonneg=True, name='l')
g = cp.Parameter(nonneg=True, name='g')
prob = cp.Problem(cp.Minimize(cp.sum_squares(X @ beta - y)     + l * cp.sum_squares(beta) + g * cp.norm(beta, 1)))

# generate code
cpg.generate_code(prob, gradient=True)

# use CVXPYlayers interface
from cvxpylayers.torch import Cvxpylayer
from cpg_code.cpg_solver import forward, backward
layer = Cvxpylayer(prob, parameters=[X,l,g], variables=[beta], custom_method=(forward, backward))
p = Gamma(w, mu, nu)  # involves beta_solution = layer(...)
p.backward()
print(w.grad, mu.grad, nu.grad)
\end{lstlisting}
\caption{Code generation and CVXPYlayers interface for elastic net example.
The integers \texttt{m}, \texttt{n}, and \texttt{J}, the constant \texttt{y},
the function \texttt{Gamma}, and the \texttt{torch} tensors \texttt{w}, \texttt{mu},
and \texttt{nu} are pre-defined.}
\label{code:elastic_CVXPYgen}
\end{figure}
The problem is modeled with CVXPY in lines~5--9.
Code is generated with CVXPYgen in line~12, where we use the \verb|gradient=True|
option to generate code for computing gradients through the problem.
After importing the CVXPYlayers interface in line 16, it is passed to the \verb|Cvxpylayer|
constructor in line 17 through the \verb|custom_method| keyword.
The performance objective is computed in line 18 and differentiated in line 19.

\paragraph{Data generation.}

We take $m = 100$ data points, $n = 20$ features,
and $J = 10$ CV folds. For every fold, we reserve $m / J = 10$ data points
for validation and use the other $90$ data points for training.
We generate the features $\bar z_i$ without outliers
by sampling from the Gaussian $\mathcal{N}(0, 1)$.
Then, we sample $\bar{\beta} \sim \mathcal{N}(0, I)$ and
set noisy labels $y_i = \bar z_i^T \bar \beta + \xi_i$ with
$\xi_i \sim \mathcal{N}(0, 0.01)$.
Afterwards, we simulate feature outliers
due to, \eg, data capturing errors.
For every feature, we randomly select $m/10 = 10$ indices
and increase the magnitude of the respective entries of $\bar z_i$ 
to a value $\sim \mathcal{U}[2, 4]$, \ie, between 2 and 4 standard
deviations, and save them in $z_i$.
We estimate the optimal cross-validated RMSE as $\hat p = 0.1$
corresponding to the standard deviation of the noise $\xi$ (if it was known).
This is a very optimistic estimate, since it implies that the data
is outlier-free after winsorization.
The tuning parameters are initialized to
$\omega^0 = (w, \mu, \nu)^0 = (3 \cdot \ones, 0, 0)$.
We set the termination tolerances to
$\epsilon^\text{rel} = \epsilon^\text{abs} = 10^{-3}$.

\paragraph{Results.}

\begin{figure}
\centering
\includegraphics[width=0.8\columnwidth]{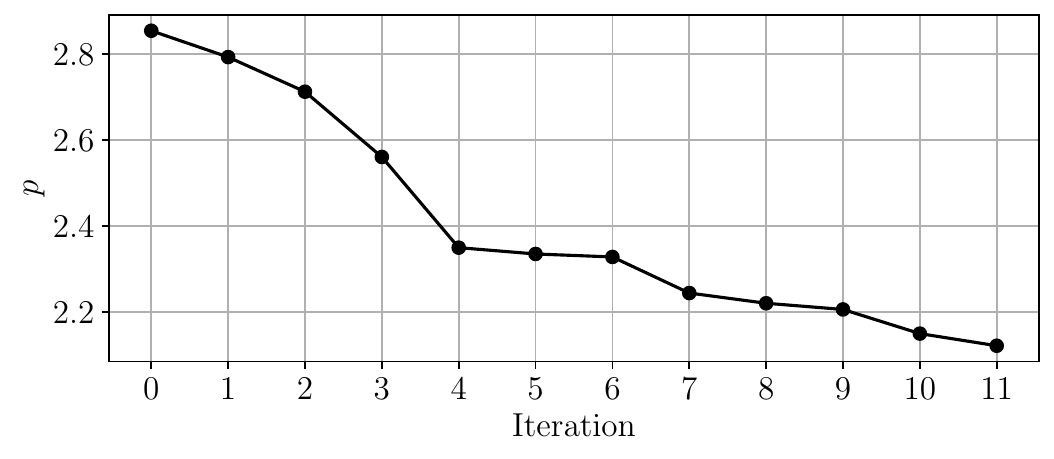}
\caption{CV loss over tuning iterations.}
\label{fig:elastic_performance}
\end{figure}

\begin{figure}
\centering
\includegraphics[width=0.8\columnwidth]{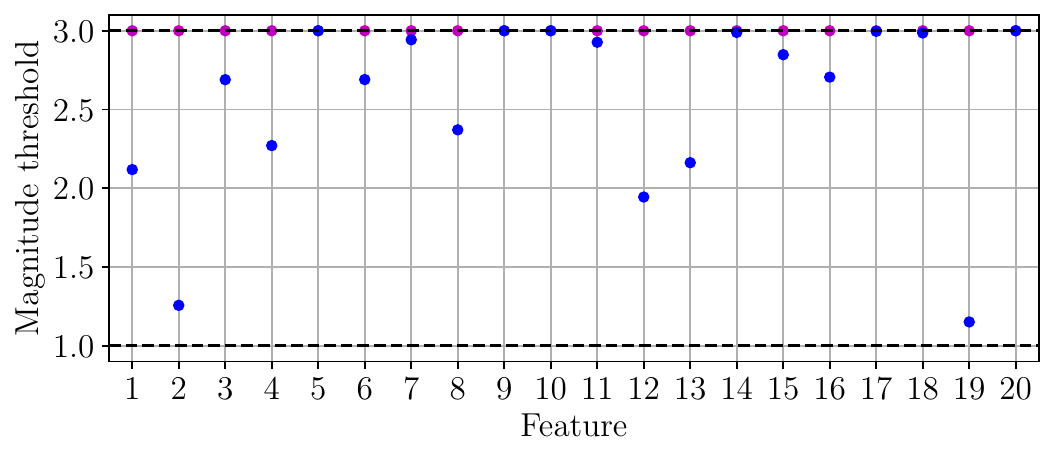}
\caption{Thresholds before (magenta) and after tuning (blue). The dashed black
lines show the range of possible winsorization thresholds.}
\label{fig:elastic_thresholds}
\end{figure}

The projected gradient method terminates after $11$ steps with a reduction
of the cross-validated RMSE from $2.85$ to $2.12$, as shown in
figure~\ref{fig:elastic_performance}.
Figure~\ref{fig:elastic_thresholds} shows the tuned winsorization thresholds $w$
and we obtain $\lambda \approx 0.68$ and $\gamma \approx 0.80$.

\paragraph{Timing.}
Table~\ref{tab:elastic_times} shows that the speed-up factor
for the gradient computations is about 5.
The speed-up for the full tuning loop is reduced due to Python overhead.
\begin{table}[h]
\begin{minipage}{\textwidth}
\centering
\begin{tabular}{lrrr}
\hline
& Full tuning & Solve and Gradient & Gradient \\ \hline
CVXPYlayers & 1.684 sec & 1.570 sec & 0.033 sec \\
CVXPYgen & 0.696 sec & 0.490 sec & 0.007 sec \\
\hline
\end{tabular}
\caption{Computation times with CVXPY and CVXPYgen for the elastic net example.}
\label{tab:elastic_times}
\end{minipage}
\end{table}

\FloatBarrier

\subsection{Approximate dynamic programming controller}\label{sec:example_ADP}

We investigate controller design with an ADP controller~\cite{wang2015approximate,
keshavarz2014quadratic} where the state and input cost parameters are subject to tuning.

We are given the discrete-time dynamical system
\BEQ\label{eq:linsys}
x_{t+1} = A x_t + B u_t + w_t, \quad t=0,1, \ldots,
\EEQ
with state $x_t \in \reals^n$, input $u \in \reals^m$
limited as $\norm{u}_\infty \leq 1$,
and state disturbance $w_t$, where $w_t$ are unknown, but assumed 
IID $\mathcal N(0,W)$, with $W$ known.
The matrices
$A \in \reals^{n \times n}$ and $B \in \reals^{n \times m}$
are the given state transition and input matrices, respectively.

We seek a state feedback controller $u_t = \phi(x_t)$ that guides the
state $x_t$ to zero while respecting the constraint $\norm{u}_\infty \leq 1$.
We judge a controller $\phi$ by the metric
\[
J = \lim_{T\to \infty} \frac{1}{T} \Expect \sum_{t=0}^{T-1} \left( 
x_t^T Q x_t + u_t^T R u_t\right),
\]
where $Q \in \symm^n_+$ and $R \in \symm^m_{++}$ are given.
We assume that the matrix $A$ contains no unstable
eigenvalues with magnitude beyond $1$, such that $J$ is guaranteed to exist.
Corresponding to \S\ref{sec:cocptuning}, we will take our performance metric as
\[
p = \Gamma(\omega) = \ell^\text{cl} (x_0, \omega)
= \sum_{t=0}^{T-1} \left( x_t^T Q x_t + u_t^T R u_t\right),
\]
where $T$ is large and fixed, and $w_t$ are sampled from $\mathcal N(0,W)$.
Note that all $x_1, \ldots, x_{T-1}$ and $u_0, \ldots, u_{T-1}$ are fully determined
by the initial state $x_0$, the controller $\phi(x; \omega)$,
and the system dynamics~\eqref{eq:linsys}.

When the input constraint is absent, we can find the optimal controller
(\ie, the one that minimizes $J$)
using dynamic programming, by minimizing a convex quadratic function,
\[
(A x_t + B u)^T P^\text{lqr} (A x_t + B u) + u^T R u,
\]
where the matrix $P^\text{lqr} \in \symm^n_{++}$ is the solution of the 
algebraic Riccati equation (ARE) for discrete time systems.
The minimizer is readily obtained
analytically, with $u$ a linear function of the state $x_t$.
See, \eg,~\cite{kwakernaak1972linear, pappas1980numerical, 	anderson2007optimal}.

We will use 
an approximate dynamic programming (ADP) controller
\[
\phi(x_t; \omega)
= \argmin_{u \in \mathcal{U}} \ell(x_t, u, \omega)
= \argmin_{\|u\|_\infty \leq 1}
(A x_t + B u)^T (P^\text{lqr} + Z) (A x_t + B u) + u^T R u.
\]
The controller is designed by $\omega = Z$ with $\Omega = \symm_+^n$.
The state $x_t$ is another parameter and the matrices
$A$, $B$, $P^\text{lqr}$, and $R$ are constants.

The quadratic form in the objective makes the above formulation non-DPP.
We render the problem DPP as
\[
\phi(x_t; \omega) = \argmin_{\|u\|_\infty \leq 1} \norm{g + H u}_2^2 + u^T R u,
\]
where the DPP parameter $\theta$ consists of
${g = L^TAx_t}$ and ${H = L^TB}$. Here, $L = \Chol(P^\text{lqr} + Z)$, where
$\Chol(\cdot)$ returns the lower Cholesky factor of its argument.
In other words, $LL^T = P^\text{lqr} + Z$ and $L$ is lower triangular.
(When running our projected gradient descent algorithm, we modify $L$ instead
of $Z$ and recover $Z = LL^T - P^\text{lqr}$ at the end of the tuning).

\paragraph{Data generation.}

We choose $n = 6$ states and $m=3$ inputs.
We consider an open-loop system $A = \diag a$ with stable and unstable modes
sampled from $[0.99, 1.00]$.
The entries of the input matrix $B$ are sampled from $[-0.01, 0.01]$.
We initialize the state at the origin and
simulate for $T=1000$ steps with noise covariance $W = 0.1^2 I$.
The \emph{true} state and input cost matrices are $Q = R = I$,
respectively, which we use to compute $P^\text{lqr}$ as the solution to the ARE.
We initialize $\omega^0 = Z^0 = 0$ and estimate the optimal control performance
$\hat p$ by running the simulation with the input constraint of the controller removed.
We use $\epsilon^\text{rel} = \epsilon^\text{abs} = 0.005$.

\paragraph{Results.}

\begin{figure}
\centering
\includegraphics[width=0.8\columnwidth]{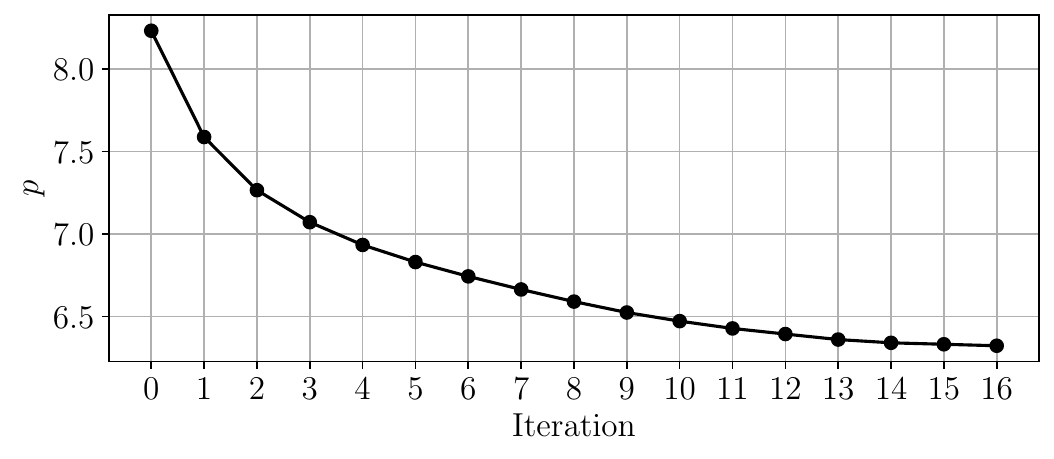}
\caption{Control performance over tuning iterations.}
\label{fig:ADP_performance}
\end{figure}

The projected gradient descent algorithm terminates after 16 gradient steps
with a reduction of the control objective from $8.23$ to $6.32$, as shown in
figure~\ref{fig:ADP_performance}.

\paragraph{Timing.}
The gradient computations are sped up compared to CVXPY
by a factor of about 40.
Including Python overhead, the whole
tuning loop is still sped up by a factor of about 6, as shown in
table~\ref{tab:ADP_times}.
\begin{table}[h]
\begin{minipage}{\textwidth}
\centering
\begin{tabular}{lrrr}
\hline
& Full tuning & Solve and gradient & Gradient \\ \hline
CVXPYlayers & 107.2 sec & 101.9 sec & 15.7 sec\\
CVXPYgen & 18.3 sec & 10.5 sec & 0.4 sec  \\
\hline
\end{tabular}
\caption{Computation times with CVXPY and CVXPYgen for the ADP controller tuning example.}
\label{tab:ADP_times}
\end{minipage}
\end{table}

\FloatBarrier

\subsection{Portfolio optimization}\label{sec:example_portfolio}

We consider a variant of the classical Markowitz portfolio optimization model
\cite{markowitz1952portfolio}
with holding cost for short positions, transaction cost, and a leverage limit
\cite{boyd2024markowitz,lobo2007portfolio}, embedded in a multi-period
trading system~\cite{boyd2017multi}.

We want to find a fully invested portfolio of holdings in $N$ assets.
The holdings are represented relative to the total portfolio value,
in terms of weights $w \in \reals^N$ with ${\ones^T w = 1}$.
The expected portfolio return is $\mu ^T w$, with esimated returns
$\mu \in \reals^N$.
The variance or risk of the portfolio return is $w^T \Sigma w$, with
estimated asset return covariance $\Sigma \in \symm^N_{++}$.
We assume that $\mu$ and $\Sigma$ are pre-computed,
which we will detail later.
We approximate the cost of holding short positions as $\kappa^\text{hold} \ones^T w_-$,
where $\kappa^\text{hold}$ describes the (equal) cost of holding a short position in 
any asset.
Subscript ``$-$'' denotes the negative part, \ie, $w_- = \max\{-w,0\}$.
We approximate the transaction cost as ${\kappa^\text{tc} \norm{w - w^\text{pre}}_1}$,
where  $\kappa^\text{tc}$ describes the cost of trading any asset
and $w^\text{pre}$ is the pre-trade portfolio.
We solve
\BEQ\label{eq:portfolio}
\begin{array}{ll}
\mbox{maximize} &\mu^T w - \gamma^\text{risk} w^T \Sigma w
-  \gamma^\text{hold} \kappa^\text{hold} \ones^T w_-
- \gamma^\text{tc} \kappa^\text{tc} \norm{\Delta w}_1\\
\mbox{subject to} &\ones^T w = 1, \quad \norm{w}_1 \leq L, \quad
\Delta w = w - w^\text{pre},
\end{array}
\EEQ
where $w, \Delta w \in \reals^{N}$ are the variables.
We introduced the variable $\Delta w$, also referred to as the \emph{trade vector},
to prevent products of parameters and render the problem DPP.
Problem~\eqref{eq:portfolio} can also be seen as a convex optimization control policy
as described in \S\ref{sec:cocptuning}, where the expected returns $\mu$
(updated once per trading period) and the previous portfolio $w^\text{pre}$ are the state $x$
and the trade $\Delta w$ is the input $u$.
Our design consists of the leverage limit $L$ and the
aversion factors $\gamma^\text{risk}, \gamma^\text{hold}, \gamma^\text{tc} > 0$
for risk, holding cost, and short-selling cost, respectively.
Since the model performance depends primarily on the orders of magnitude of
these factors, we write them as
\[
\gamma^\text{risk} = 10^{\nu^\text{risk}}, \quad
\gamma^\text{hold} = 10^{\nu^\text{hold}}, \quad
\gamma^\text{tc} = 10^{\nu^\text{tc}},
\]
and tune $\omega = (L, \nu^\text{risk}, \nu^\text{hold}, \nu^\text{tc})$,
restricted to the design space
\[
\Omega = [1, 2] \times [-3, 3]^3. 
\]
We keep the risk $\Sigma$ and costs $\kappa^\text{hold}$ and $\kappa^\text{tc}$
constant.

We evaluate the performance of the model via a back-test over $h$ trading periods.
After solving problem \eqref{eq:portfolio} at a given period,
we trade to $w^\star$, pay short-selling cost
$\kappa^\text{hold} \ones^T w^\star_-$ and transaction cost
$\kappa^\text{tc} \norm{w^\star - w^\text{pre}}_1$, experience the returns $r_t$,
and re-invest the full portfolio value.
Hence, the total portfolio value evolves as
\[
V_{t+1} = V_t (1 + r_t^T w^\star) - \kappa^\text{hold} \ones^T w^\star_-
- \kappa^\text{tc} \norm{w^\star - w^\text{pre}}_1.
\]
The pre-trade portfolio for the following trading period is re-balanced as
\[
w^\text{pre} = w^\star \circ (1 +  r_t) \cdot V_t / V_{t+1}
\]
and the portfolio realized return at period $t$ is
\[
R_t = (V_{t+1}- V_t) / V_t.
\]
We consider the average return and portfolio risk,
\[
\bar R = (1/h) \sum_{t = 1}^h R_t, \quad
\sigma =  \left( (1/h) \sum_{t=1}^h R_t^2\right)^{1/2},
\]
respectively, and annualize them as
\[
\bar R^\text{ann} = h^\text{ann} \bar R, \quad
\sigma^\text{ann} = (h^\text{ann})^{1/2} \sigma,
\]
where $h^\text{ann}$ is the number of trading periods per year.
We take their ratio as performance metric, the so-called
\emph{Sharpe ratio (SR)}~\cite{ledoit2008robust},
\[
p = \text{SR} = \bar R^\text{ann} / \sigma^\text{ann}
= (h^\text{ann})^{1/2} \bar R / \sigma.
\]

\paragraph{Data generation.}

We consider three adjacent intervals of trading periods. First, we use a \emph{burn-in}
interval to compute the estimate for the expected returns at later time periods
and to compute the constant risk estimate. Second we take a \emph{tune} interval
to perform parameter optimization. Third, we use a \emph{test} interval
to evaluate the final parameter choice out-of-sample.
We denote the lengths of the three intervals by $h^\text{burnin}$, $h^\text{tune}$,
and $h^\text{test}$, respectively.

We compute the expected return $\mu_t$
for the holdings at time $t$ as the back-looking moving average of
historical returns $r_t$ with window size $h^\text{burnin}$.
To compute the constant risk estimate, we first compute the empirical
covariance $\hat\Sigma$ of returns over the burn-in interval.
Then, we fit the standard factor model
\[
\Sigma = FF^T+D
\]
to $\hat \Sigma$, where $F \in \reals^{N \times K}$ is the factor loading matrix and
the diagonal matrix $D \in \symm^N_{++}$
stores the variance of the idiosyncratic returns~\cite{elton2009modern}.

We consider $N=25$ stock assets, chosen randomly from the S\&P 500,
where historical return data is available from 2016--2019.
While this clearly imposes survivership bias~\cite{brown1992survivorship},
the point of this experiment is not to find realistic portfolios
but rather to assess the numerical performance of CVXPYgen.
We choose $K=5$ factors,
$h^\text{burnin} = 260$, $h^\text{tune} = 520$, 
$h^\text{test} = 260$, and we take the number of trading periods per year
as $h^\text{ann} = 260$, \ie, we trade once a day.
In other words, we use data from the year 2016 as burn-in interval
for estimating $\mu_t$ and to estimate $\Sigma$.
We tune the model with data
from the years 2017 and 2018 and test it with data from the year 2019.
We fix ${\kappa^\text{hold} = \kappa^\text{tc} = 0.001}$.
We estimate the optimal Sharpe ratio to be roughly $\hat p = 1$.
We initialize the design vector as
$\omega^0=
(L, \nu^\text{risk}, \nu^\text{hold}, \nu^\text{tc})^0 = (1, 0, 0, 0)$.
We set the termination tolerances to
$\epsilon^\text{rel} = \epsilon^\text{abs} = 0.03$.

\paragraph{Results.}

\begin{figure}
\centering
\includegraphics[width=0.8\columnwidth]{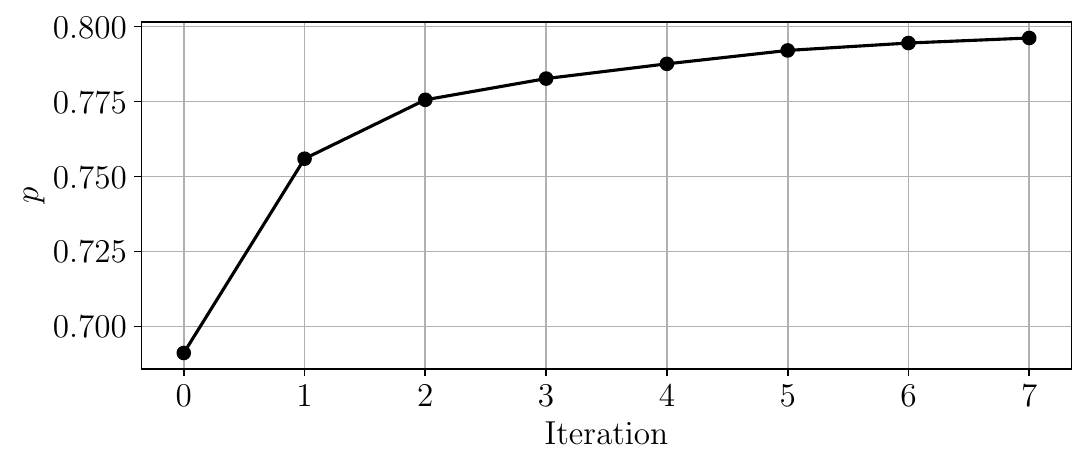}
\caption{Sharpe ratio over tuning iterations.}
\label{fig:portfolio_performance}
\end{figure}

The projected gradient descent algorithm terminates after 7 iterations
and improves the Sharpe ratio by a bit more than $0.1$, from $0.69$ to $0.80$,
where it is saturating, as shown in figure~\ref{fig:portfolio_performance}.

\begin{figure}
\centering
\includegraphics[width=0.8\columnwidth]{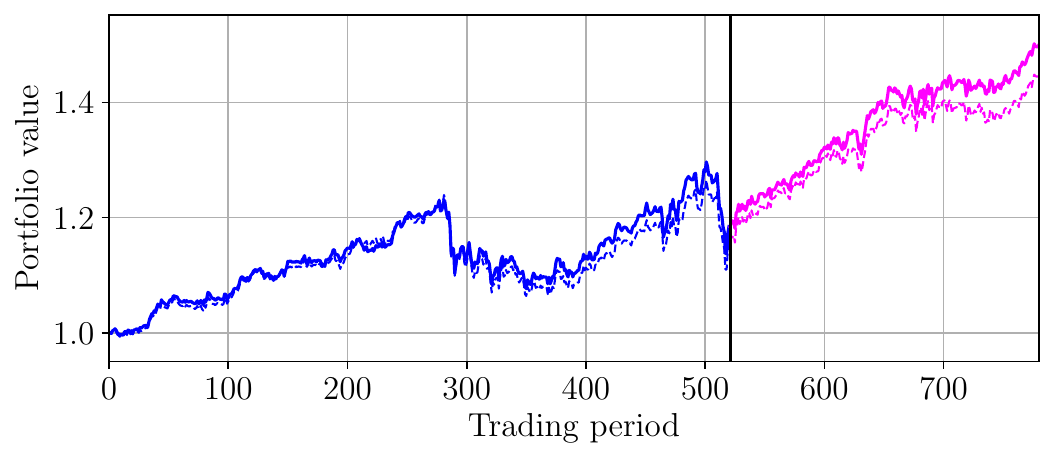}
\caption{Portfolio value evolution before (dashed line) and after tuning (solid line).
Blue and pink color represent the tuning and testing intervals, respectively.}
\label{fig:portfolio_trajectories}
\end{figure}

The values of the tuning parameters are changed to
$\gamma^\text{risk} \approx 10$, $\gamma^\text{hold} \approx 1$,
$\gamma^\text{tc} \approx 1.2$, and $L \approx 1$.
Figure~\ref{fig:portfolio_trajectories} contains the portfolio value
over the trading periods used for tuning and out-of-sample,
before and after tuning, respectively.

Table~\ref{tab:sharpe_ratios} contains the respective Sharpe ratios.
While the tuning interval appears to be a difficult time period with
large drawdown in the middle and the end of the interval,
the Sharpe ratio is improved out-of-sample from an already high level.

\begin{table}[h]
\begin{minipage}{\textwidth}
\centering
\begin{tabular}{lrr}
\hline
& In-sample & Out-of-sample \\ \hline
Before tuning & 0.69 & 2.21  \\
After tuning & 0.80 & 2.39 \\
\hline
\end{tabular}
\caption{Sharpe ratios.}
\label{tab:sharpe_ratios}
\end{minipage}
\end{table}

\FloatBarrier

\paragraph{Timing.}
Table~\ref{tab:portfolio_times} shows the solve and differentiation times.
The gradient computations are sped up by a factor of about 10.
Including Python overhead, the overall tuning loop is sped up by
a factor of about 3.
\begin{table}[h]
\begin{minipage}{\textwidth}
\centering
\begin{tabular}{lrrr}
\hline
& Full tuning & Solve and gradient & Gradient \\ \hline
CVXPYlayers & 61 sec & 57 sec & 23 sec \\
CVXPYgen & 21 sec & 17 sec & 2 sec \\
\hline
\end{tabular}
\caption{Computation times with CVXPY and CVXPYgen for the portfolio optimization example.}
\label{tab:portfolio_times}
\end{minipage}
\end{table}

\FloatBarrier

\section{Conclusions}\label{sec:conclusion}

We have added new functionality to the code generator CVXPYgen
for differentiating through parametrized convex optimization problems.
Users can model their problem in CVXPY with instructions close to the math, and create an
efficient implementation of the gradient computation in C,
by simply setting an additional keyword argument of the \mbox{CVXPYgen}
code generation method.
Our numerical experiments show that the gradient computations
are sped up by around one order of magnitude for typical use cases.

{\small
	\bibliography{refs}
}

\end{document}